\documentclass[12pt]{amsart}

\usepackage{epsf,amssymb}

\newtheorem{thm}{Theorem}[section]

\newtheorem*{claim}{Claim}

\newtheorem{prop}[thm]{Proposition}
\newtheorem{defn}[thm]{Definition}
\newtheorem{lemma}[thm]{Lemma}

\newtheorem{conj}[thm]{Conjecture}

\newcommand{\R}{\mathbb{R}}
\newcommand{\Z}{\mathbb{Z}}

\renewcommand{\P}{\mathbb{P}}
\newcommand{\bdry}{\partial}

\newcommand{\sa}{\rightsquigarrow}
\renewcommand{\sc}{\rightarrowtail}
\newcommand{\s}{\vskip.1in}
\newcommand{\n}{\noindent}

\newcommand{\be}{\begin{enumerate}}
\newcommand{\ee}{\end{enumerate}}

\begin{document}
\title{Convex decomposition theory}

\author{Ko Honda}
\address{University of Georgia, Athens, GA 30602 and AIM, Palo Alto, CA 94306} 
\email{honda@math.uga.edu}
\urladdr{http://www.math.uga.edu/\char126 honda}

\author{William H. Kazez}
\address{University of Georgia, Athens, GA 30602}
\email{will@math.uga.edu}
\urladdr{http://www.math.uga.edu/\char126 will}

\author{Gordana Mati\'c}
\address{University of Georgia, Athens, GA 30602}
\email{gordana@math.uga.edu}

\date{This version: January 31, 2001.}

\thanks{KH supported by NSF grant DMS-0072853 and the American Institute of Mathematics;
GM supported by NSF grant DMS-0072853; WHK supported by NSF grant DMS-0073029.}

\begin{abstract}
	We use convex decomposition theory to (1) reprove the existence of a universally tight
	contact structure on every irreducible 3-manifold with nonempty boundary, and
	(2) prove that every toroidal 3-manifold carries infinitely many nonisotopic, nonisomorphic
	tight contact structures.
\end{abstract}
\maketitle

It has been known for some time that there are deep connections
between the theory of taut foliations
and tight contact structures due to the work of Eliashberg and
Thurston \cite{ET}.  In particular, they proved
that a taut foliation can be perturbed into a (universally) tight
contact structure.
In a previous paper \cite{HKM} we introduced the notion of {\it
convex decompositions} and explained
how convex decompositions can naturally be viewed as generalizations
of {\it sutured manifold decompositions}
introduced by Gabai in \cite{Ga}.   In this paper we take the
viewpoint that convex decompositions
are completely natural in 3-dimensional contact topology and that
many theorems can be proven
directly in the category of tight contact manifolds with convex
splittings as morphisms.

The first theorem of the paper is a version of a theorem
by Gabai-Eliashberg-Thurston in the case of manifold with boundary:

\begin{thm} [Gabai-Eliashberg-Thurston]\label{thm1} Let $(M,\gamma)$ be an oriented, compact,
connected, irreducible, sutured 3-manifold which has boundary, is taut, and has annular sutures. 
Then $(M,\gamma)$ carries a universally tight contact structure.
\end{thm}

We provide an alternate proof which (1) does not require us to perturb taut foliations into tight
contact structures  as we did in \cite{HKM} and (2) does not  resort to four-dimensional {\it
symplectic filling} techniques in order to prove tightness.  Instead, we use Gabai's sutured
manifold decomposition and directly apply a gluing theorem (Theorem \ref{gluing}). To prove the
gluing theorem we apply key ideas from Colin's papers \cite{Co99a} and \cite{Co99b} in the
context of convex decomposition theory.

We also apply similar ideas to prove the following theorem:

\begin{thm} \label{infinite}
Let $M$ be an oriented, closed, connected, irreducible 3-manifold which contains
an incompressible torus.
Then $M$ carries infinitely many isomorphism classes of universally
tight contact structures.
\end{thm}

This theorem confirms a conjecture which has its beginnings in the works of
Giroux \cite{Gi94} and Kanda \cite{K97}, was proved for torus bundles
over $S^1$ by Giroux \cite{Gi99a}, and
was extended by Colin \cite {Co99b} to the case where there exist two
incompressible tori $T$, $T'$ with
a ``persistent'' intersection. The flip side of Theorem \ref{infinite}
would be the following conjecture:

\begin{conj}
Let $M$ be a closed, connected, irreducible 3-manifold which is
atoroidal (does not contain an
incompressible torus).  Then $M$ carries only finitely many isotopy
classes of tight contact structures.
\end{conj}

It should be noted that this work has been influenced greatly by the
work of Colin \cite{Co99a, Co99b,
Co99c}. In particular, the reader familiar with Colin's work will
recognize that many of the ideas of
this paper are adaptations and strengthenings of Colin's ideas in the
setting of convex decompositions.   We have also
been informed by Colin that he has recently obtained Theorem
\ref{infinite} (independently).

\s
\n
{\it Necessary background.} This paper is intended as a sequel to
\cite{HKM}, and it is recommended
that the reader read it first.  We will freely use the circle of
ideas introduced there.  The reader will
also find it helpful to have read \cite{H1} and \cite{H2} and to have
familiarized him/herself with ideas
of gluing and bypasses.

\s
\n
{\it Conventions.}
\be
\item $M$ = oriented, compact $3$-manifold.
\item $\xi$ = positive contact structure which is co-oriented by a
global 1-form $\alpha$.
\item A convex surface $\Sigma$  is either closed or compact with
Legendrian boundary.
\item $\Gamma_\Sigma$ = dividing set of a convex surface $\Sigma$.
\item $\#\Gamma_\Sigma$ = number of connected components of $\Gamma_\Sigma$.
\item $|\beta\cap \gamma|$ = geometric intersection number of two
curves $\beta$ and $\gamma$
on a surface.
\item $\#(\beta\cap\gamma)$ = cardinality of the intersection.
\item $\Sigma\backslash \gamma$ = metric closure of the complement of
$\gamma$ in $\Sigma$.
\item $t(\beta,Fr_S)$ = twisting number of a Legendrian curve with respect to the framing 
induced from the surface $S$.
\ee

\section{Gluing and proof of  Theorem \ref{thm1}}

In this paper we will assume that the 
reader is familiar with concepts of sutured manifold decompositions and convex decompositions. 
In particular, we refer the reader to \cite{HKM} for an explanation of the relationship between 
sutured manifold decompositions and convex decompositions. We will
also need the following definitions.

\begin{defn}  A {\rm sutured manifold with annular sutures}
is a sutured manifold
$(M, \gamma)$ which satisfies the following:
\be
\item $\bdry M$ is nonempty.
\item Every component of $\bdry M$ contains an annular suture.
\item Every component of $\gamma$ is an annulus, that is $\gamma =
A(\gamma) \ne \emptyset$.
\ee
\end{defn}

\begin{defn}  Let $S$ be a properly embedded compact convex surface
with Legendrian boundary.
A connected component $\gamma$ of $\Gamma_S$ is {\em
boundary-parallel} (or $\bdry$-parallel) if it  is an arc with boundary on $\bdry S$ that
cuts off a half-disk $D$ from $S$.  We also
say that $\Gamma_S$ is
{\em $\bdry$-parallel} if all of its connected components are
$\bdry$-parallel.
\end{defn}

Let $(M_0,\gamma_0)$ be a sutured manifold with annular sutures and let
$$(M_0,\gamma_0)\stackrel{S_0}{\sa} (M_1,\gamma_1)\stackrel{S_1}{\sa}\cdots
\stackrel{S_{n-1}}{\sa}(M_n,\gamma_n)=\cup (B^3,S^1 \times I),$$ be a
sutured manifold
hierarchy.  We define the corresponding convex hierarchy
$$(M_0,\Gamma_0)\stackrel{(S_0,\sigma_0)}{\sa}
(M_1,\Gamma_1)\stackrel{(S_1,\sigma_1)}{\sa}\cdots
\stackrel{(S_{n-1},\sigma_{n-1})}{\sa}(M_n,\Gamma_n)=\cup
(B^3,S^1),$$ as follows.  Let
$(M_i, \Gamma_i)$ be the convex structure associated to
$(M_i,\gamma_i)$, that is $\Gamma_i$
consists of the oriented cores of $\gamma_i$.  We may assume that
each component of $S_i$
which is contained in $\gamma_i$ intersects $\Gamma_i$ transversely
in at least two points.
Each $\sigma_i$ is defined to be a collection of {\it $\bdry$-parallel}
arcs in $S_i$ such that
the half-disks they bound in $S_i$ intersect $\Gamma_i$ in one point each
and have the opposite
orientation as $S_i$, and the orientation induced by the half-disks
on the components of
$\sigma_i$ cause $\sigma_i$ to start in $R_-(\Gamma_i)$ and end in $R_+(\Gamma_i)$.  With
these conventions it is
straightforward to see that $(M_i, \Gamma_i)$ split along $(S_i, \sigma_i)$ is
$(M_{i+1},\Gamma_{i+1})$.

\begin{figure} 	{\epsfysize=2in\centerline{\epsfbox{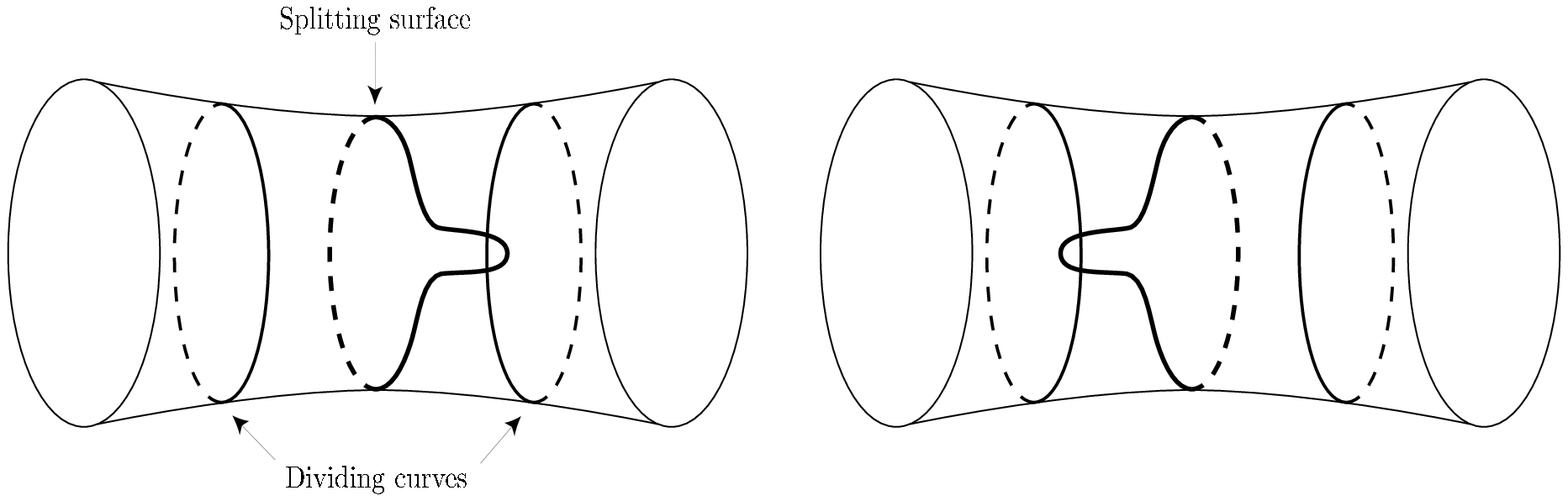}}}
	\caption{} 	\label{finger}
\end{figure}

For the purposes of a later argument we need to know that we can assume that
each component of each
$\partial S_i$  intersects
the dividing curves $\Gamma_i$ at least twice.  To show how this may be
arranged, consider the
important special case where some component $R$ of
$\bdry M_i \backslash \Gamma_i$ contains
a component $C$ of $\partial S_i$ such that $C$ does not separate
$R$, and $R$ intersects no
other components of $\partial S_i$.  $S_i$ may be perturbed in a
neighborhood of $C$ by
a ``finger move'' in two different ways, ${S_i}'$, ${S_i}''$ (see
Figure \ref{finger}).  As
described above, add appropriate $\bdry$-parallel dividing curves to $S_i'$ and
$S_i''$.  Splitting along
precisely one of $(S_i',\sigma_i')$ and $(S_i'',\sigma_i'')$ will
produce the convex structure
$(M_{i+1},\Gamma_{i+1})$, and we continue to abuse notation by
denoting this splitting
surface $(S_i, \sigma_i)$.  This technique may also be applied if $C$
is replaced by a family
of parallel coherently oriented curves in $R$.  In general, by
Theorem \ref{Gabai:hierarchy} stated below,
we may always assume this is the case by choosing the original
sequence $S_1, S_2, \dots$ to
be``well-groomed''.

\begin{thm}[Gabai] \label{Gabai:hierarchy} Let $(M,\gamma)$ be a taut
sutured manifold with
$H_2(M,\bdry M)\not =0$. Then $(M,\gamma)$ admits a sutured manifold
hierarchy such that
$S_i\cap \bdry M_{i}\not=\emptyset$ if $\bdry M_i\not=\emptyset$, and
$S_i$ is well-groomed
along $\bdry M_i$, i.e., for every component $R\subset \bdry
M_i\backslash \gamma_i$,
$S_i\cap R$  is a union of parallel oriented nonseparating
simple closed curves if
$R$ is nonplanar and arcs if $R$ is planar.
\end{thm}

We then have the following theorem:

\begin{thm}    \label{glueback} Let $(M,\gamma)$ be a sutured
manifold with annular sutures.
If $(M,\gamma)$ admits a sutured manifold hierarchy, then there
exists a universally
tight contact structure on $M$ with convex boundary and dividing set
$\Gamma$ corresponding
to $\gamma$.
\end{thm}

\proof
  Note that existence of a sutured manifold hierarchy implies that 
$(M,\gamma)$ is taut.
Use Theorem
\ref{Gabai:hierarchy} to choose a well-groomed sutured  manifold 
hierarchy for $(M,
\gamma)$ and consider the associated convex hierarchy of $(M,
\Gamma)$.  Starting with
the last term of the convex hierarchy, $\cup(B^3,S^1)$,
inductively define contact
structures on the terms of the convex hierarchy. For the initial
step, we use the
following fundamental result of Eliashberg
\cite{E92}:

\begin{thm}[Eliashberg]
\label{thm:unique} Assume there exists a contact structure $\xi$ on a
neighborhood of
$\bdry B^3$ which makes $\bdry B^3$ convex with $\#\Gamma_{\bdry B^3}=1$.
Then there exists a unique extension of $\xi$ to a tight contact
structure on $B^3$
up to an isotopy which fixes the boundary.
\end{thm}

Now that we have unique tight contact structures on the 3-balls, we
use the Gluing Theorem (Theorem \ref{gluing}) below to
successively obtain a universally tight contact structure on the
glued-up manifold.  The
hypotheses that $(M,\gamma)$ is a sutured manifold with annular
sutures and the choice of
perturbation of the well-groomed splitting surfaces used in the 
convex hierarchy of
$(M,\Gamma)$ guarantee that the hypotheses of Theorem \ref{gluing} are
satisfied. \qed

Note that a splitting surface of a convex splitting associated to a sutured manifold 
splitting has a $\bdry$-parallel dividing set.  This is crucial as we are about to use
the following Gluing Theorem to finish our proof. This theorem first 
appeared in  Colin
\cite{Co99a} in a slightly different form. Our formulation is in 
terms of convex
decompositions, and our proof  relies on bypasses and dividing curve 
configurations.

\begin{thm} [Gluing]  \label{gluing}  Let $(M,\xi)$ be an irreducible contact
manifold with nonempty convex boundary, and let $S\subset M$ be a 
properly embedded
compact convex surface with Legendrian boundary such that (1) $S$ is 
incompressible
in $M$, (2) $t(\gamma,Fr_S)<0$ for each connected component $\gamma$ of
$\bdry S$, i.e., each component of $\gamma$ nontrivially intersects the dividing set
$\Gamma_{\bdry M}$, and (3) $\Gamma_S$ is $\bdry$-parallel. Consider a
decomposition  of $(M,\xi)$ along $S$. If $\xi$ is universally tight on
$M\backslash S$,  then $\xi$ is universally tight on $M$.
\end{thm}

\begin{proof}
We will prove that, under conditions (1)--(3), the fact that $\xi$ is universally tight on
$M\backslash S$ implies that $\xi$ is tight on $M$.  The proof, as will be
clear, easily adapts to any finite cover of $M$, hence showing that 
the universal
cover is tight, provided $M$ is {\it residually finite} (which is the case, since $M$ is 
Haken --- see Section \ref{section:residual} for a discussion of residual finiteness). Assume, 
on the contrary, that $M$ is not tight and that $D$ is an overtwisted disk in $M$. Then, after 
a possible contact isotopy, we can assume that $D$ intersects $S$ transversally along 
Legendrian curves and arcs  and that $\bdry D \cap S \subset
\Gamma_S$  (see Lemma 2.7 in \cite{H2}). Note that closed curves in $D\cap S$
are homotopically trivial on $S$ since $S$ is incompressible.  We would like to
argue that, starting with innermost closed curves, we can eliminate them  by
pushing $S$ across $D$ (we do not push $D$ across $S$ to avoid introducing
self-intersections into $D$). Since $M$ is assumed to be irreducible, the 2-sphere
formed by two disks (one on $D$ and one on $S$) bounding an innermost curve of intersection 
$\delta$ on $S$  bounds a ball across which $S$ can be isotoped.
We will include the sketch of the following fact (see Lemma~2.8 
and~2.9 in \cite{H2}),
as we will need to refer to some steps in the proof of it later.

\begin{lemma} We can push $S$ across $D$ to eliminate $\delta$ in a 
finite number
of steps, each of which is a bypass along an arc of the circle $\delta$.
\end{lemma}

\begin{proof} Consider the subdisk  $D_{\delta}$ of $D$ bounded by 
$\delta$. Since $\delta$ is a homotopically trivial Legendrian curve on $S$, 
$t(\gamma,Fr_S)$ must be negative.  After possible perturbation rel boundary, $D_{\delta}$ is 
convex with Legendrian boundary satisfying $t(\bdry D_{\delta})<0$.  $\Gamma_{D_{\delta}}$ 
consists of properly embedded arcs with endpoints on $\bdry D_\delta$ and  no closed  
components, since the interior of $D_{\delta}$ lies in the tight manifold $M\backslash S$. We 
push $S$ to engulf a bypass along $\delta$ corresponding to a $\bdry$-parallel dividing curve 
on $D_{\delta}$. We can continue
until there is only one arc left. It is not  hard then to see that 
the last bit of
the isotopy is a contact isotopy.
\end{proof}

In a similar manner we can eliminate outermost arcs of intersection, 
by engulfing
the half-disk of $D$ bounded by such an arc in a ball, and then 
pushing $S$ across
that  ball as in the previous lemma. Therefore, starting with innermost circles
$\delta$ and outermost arcs, we can push $D$ across $S$ to reduce 
$\#(S\cap D)$ in
steps, at the expense of  modifying the dividing curve configuration on $S$.
Changes in dividing curve configurations can be reduced to a sequence of bypass
attachment moves which must now be analyzed.

Without the assumption that the dividing set on the cutting surface consists of
$\bdry$-parallel arcs the following lemma would not be true.

\begin{lemma} Let $S$ be a convex surface with Legendrian boundary in a contact
manifold $(M,\xi)$, such that $\Gamma_S$ is $\bdry$-parallel and $(M\backslash
S, \xi)$ tight. Then any convex surface $S'$ obtained from $S$ by a 
sequence of
bypasses will have $\Gamma_{S'}$ obtained from $\Gamma_S$ by possibly adding pairs of
parallel nontrivial curves (up to isotopy rel boundary).
\end{lemma}

\begin{proof}
This follows by examing the possible bypasses.  One possibility is 
that the bypass
is trivial, that is it is contained in a
disk and produces no change in the dividing curves. Some bypasses 
will produce an
overtwisted disk,  and hence cannot exist inside a tight manifold. 
Some bypasses
introduce a pair of parallel curves, and finally a bypass may change 
one pair of
parallel curves into another pair or remove a pair of parallel curves.
\end{proof}

We need one more simple lemma:

\begin{lemma} \label{equal} Let $S$ be a convex surface with 
Legendrian boundary in
a contact manifold $(M,\xi)$, such that $\Gamma_S$ is $\bdry$-parallel. If a
convex surface $S'$ is obtained from $S$  by a bypass such that 
$\Gamma_{S'}$ is
isotopic to $\Gamma_S$, then $S$ and $S'$ are contact isotopic, and 
in particular
$(M\backslash S,\xi)$ is tight if and only if $(M\backslash S',\xi)$ is.
\end{lemma}     

\begin{proof} Left to the reader.
\end{proof}

Consider a single bypass move on $S$ with $\Gamma_S$ $\bdry$-parallel.  It is either
trivial or increases $\#\Gamma$.  If $\#\Gamma$ is increased, the arc 
of attachment
$\delta$ starts on an arc $l$ of $\Gamma_{S}$ and comes back to $l$, 
thereby generating a
non-trivial element of $\pi_1(S,l)$. There exists a large enough finite cover
$\pi:\widetilde M\rightarrow  M$ which ``expands'' $S$ to $\widetilde 
{S} =\pi^{-1}S$, such
that a lift of $\delta$ becomes a trivial arc of attachment, that is, 
it connects two
different components of $\Gamma_{\widetilde S}$.  The existence of the large finite cover 
follows from the facts that  (1) an incompressible surface $S$ has $\pi_1(S)$ which injects 
into $\pi_1(M)$ and (2) a Haken 3-manifold has residually finite $\pi_1$.

Therefore, by passing through a finite succession of covers, we  
construct a cover
$\widetilde {M}$, together with the preimage $\widetilde S$ of $S$ and a lift $\widetilde D$ of 
the overtwisted disk $D$,  in which all the bypasses needed to isotop $\widetilde S$ across 
$\widetilde D$ are trivial.  The surface $\widetilde S'$ with $\widetilde S'\cap 
\widetilde D=\emptyset$, obtained from $\widetilde S$ via trivial bypass attachments, has the 
same dividing set as $\widetilde S$ and is contact isotopic to $\widetilde S$ by Lemma 
\ref{equal}.  It follows that $(\widetilde {M} \backslash
\widetilde {S'},\pi^*\xi)$ is tight.  This contradicts the existence of an overtwisted disk 
$\widetilde D\hookrightarrow \widetilde M\backslash \widetilde S'$ and finishes the proof of 
Theorem \ref{gluing}. 
\end{proof}

Combining Theorem \ref{Gabai:hierarchy} with Theorem \ref{glueback}, we have
Theorem \ref{thm1}.

\bigskip
\section{Tight contact structures on $T^2 \times I$ and folding}

\subsection{Universally tight contact structures on $T^2\times I$}  \label{section2-1}

In this section we review key properties of universally tight contact structures on $T^2\times 
I$.  Details of assertions can be
found in \cite{H1} (also see \cite{Gi99b}).   

Fix an oriented identification 
$T^2\stackrel{\sim}{\rightarrow} \R^2/\Z^2$.  Consider $T^2\times I=T^2\times[0,1]$ with 
coordinates $((x,y),z)$.   We will focus on tight contact 
structures on $T^2 \times I$ with convex boundary which satisfy 
\begin{equation}    \label{eq1}
\mbox{slope}(\Gamma_{T_0}) = \mbox{slope}(\Gamma_{T_1}).
\end{equation}
Here we write $T_z=T^2\times \{z\}$.  The following proposition is a consequence of the 
classification of tight contact structures on $T^2\times I$.

\begin{prop} All tight contact structures satisfying Equation (\ref{eq1}) are
universally tight.
\end{prop}

The prototypical universally tight contact structures on $T^2\times I$ are $\xi_k$, $k\in 
\Z^+$, given by 1-forms $\alpha_k= \sin(\pi kz)dx+\cos(\pi kz)dy$, with the boundary adjusted 
so it becomes convex with $\#\Gamma_{T_i}=2$, $i=0,1$. If we pick a curve that meets each
component of $\Gamma_{T_0}$ in a point and take the $l$-fold cover corresponding to this
curve, we call the pull-back contact structure $\xi_k$ as well.  Note that for this version of
$\xi_k$ we have $\#\Gamma_{T_i}=2l$.

We say that a curve $C$ on a convex surface $\Sigma$ is {\it efficient} (with respect to 
$\Gamma_\Sigma$) if $|C\cap \Gamma_\Sigma|=\#(C\cap\Gamma_\Sigma)$. Given a tight contact
manifold $(T^2\times I,\xi)$ satisfying Equation (\ref{eq1}), we say an  annulus $A\subset
T^2\times I$ is {\it horizontal} if $A$ is a properly embedded, convex annulus with Legendrian
boundary, and $\bdry A$ is efficient with respect to $\Gamma_{T_0}\cup 
\Gamma_{T_1}$. Horizontal annuli encode much (sometimes all) of the information on the tight 
contact manifold $(T^2\times I,\xi)$.  

\begin{prop}    \label{horizontal}
Consider $(T^2\times I, \xi_k)$ above with $\Gamma_{T_i}=2l$, $i=0,1$. Let  $A=S^1\times 
\{pt\}\times I$ be a horizontal annulus which minimizes $\#\Gamma_A$ amongst horizontal annuli 
isotopic to $A$.  Then $\Gamma_A$ consists of $k-1$ closed curves and $l$ $\bdry$-parallel arcs 
abutting each boundary component of the annulus.  Moreover, any other horizontal annulus $A'$ 
isotopic to $A$ will have $\Gamma_{A'}$ which can be obtained from $\Gamma_A$ by adding $2m$ 
extra closed curves, $m\in \Z^{\geq 0}$.   In particular, distinct $\xi_k$'s are distinguished 
by  minimal horizontal annuli.
\end{prop}

Next, $(T^2\times I,\xi)$ satisfying Equation (\ref{eq1}) is called {\it rotative} if there is 
a convex torus $T'$ parallel to $T_0$ for which 
$\mbox{slope}(\Gamma_{T'}) \not = \mbox{slope}(\Gamma_{T_0})$.  $(T^2\times I,\xi)$ which is 
not rotative is {\it non-rotative}.  The following proposition explains 
the relationship between non-rotative $\xi$ and their horizontal annuli $A$.

\begin{prop} \label{non-rotative} If $\xi$ satisfies Equation (\ref{eq1}) and
$A$ is a horizontal annulus, then $\xi$ is non-rotative if and only if there
exists a dividing curve of $\Gamma_A$ which is {\it nonseparating} (i.e.,
connects the two boundary components of $A$).   Moreover, for a non-rotative
$\xi$, the isotopy class of the contact structure (rel boundary) is determined
completely by $\Gamma_A$.  In particular, $\xi$ is contactomorphic to an
$S^1$-invariant contact structure on $S^1\times A$, where $\{pt\}\times A$ has
dividing set $\Gamma_A$.
\end{prop}

An important consequence of Proposition \ref{non-rotative} is that
non-rotative $(T^2\times I,\xi)$ are in 1-1 correspondence with horizontal
annuli $A$ which are nonseparating.  Now, according to Giroux's
criterion \cite{Gi00a} (also see \cite{H1}), a convex surface $\not=S^2$ has a tight 
neighborhood if and only if there is no homotopically trivial dividing curve.  (Moreover, the 
tight neighborhood is automatically universally tight.)   Therefore, two non-rotative, 
universally tight $(T^2\times [0,1],\xi_1)$ and $(T^2\times[1,2],\xi_2)$ glue into a 
universally tight contact structure if and only if their corresponding horizontal annuli $A_1$ 
and $A_2$ do not glue to give a homotopically trivial dividing curve.

On the other hand, any rotative $\xi$ is isomorphic to a universally tight contact structure 
obtained by taking $(T^2\times [0,1],\xi_k)$ for some $k\in \Z^+$ and attaching non-rotative 
outer layers $T^2\times [-1,0]$ and $T^2\times[1,2]$. 

Now, consider a contact manifold $(M,\xi)$ and a torus $T\subset M$.  The {\it torsion} of the 
isotopy class $[T]$ of $T$, due to Giroux \cite{Gi99a}, is defined as follows.  Let $tor([T], 
M,\xi)$ be the supremum over $k\in \Z^{\geq 0}$, for which there exists a contact embedding 
$(T^2\times I,\xi_k) \hookrightarrow (M,\xi)$ where the image of $T^2\times\{pt\}$ is isotopic 
to $T$.  The following proposition is a consequence of the work of Giroux \cite{Gi99a} and 
Kanda \cite{K97}.

\begin{prop}  $tor([T_0],T^2\times I,\xi_k)=k$.
\end{prop}

\subsection{The Attach = Dig Principle}   \label{a=d}
Let $(T^2\times [0,1],\xi)$ be a rotative universally  tight contact structure with 
$\#\Gamma_{T_1}>>2$.   We explain the effect of attaching a bypass along
$T_1$ from the outside along a Legendrian arc $a$ which connects  3 different  components  of
$\Gamma_{T_1}$.   This operation  reduces the number of components of $\Gamma_{T_1}$  by two
and corresponds to attaching a non-rotative $(T^2\times [1,2],\xi')$  layer. If we extend the
attaching arc $a$ of the bypass to a closed Legendrian curve $\alpha$ which is efficient  with
respect to $\Gamma_{T_1}$ and intersects each component of $\Gamma_{T_1}$ exactly once,   the
dividing set that $\xi'$ induces on the horizontal annulus $ A' = \alpha \times [1,2]$ 
consists of an arc connecting  two points on $a$ lying between  the dividing curves of $T_1$,
and nonseparating arcs emanating from  all other midpoints of segments of $\alpha \backslash
\Gamma_{T_1}$.  Now we make the following simplifying assumption (which is satisfied by
$\xi_k$):

\s\n
{\bf $\bdry$-Parallel Assumption.} 
{\it The dividing set of $\xi$ on a horizontal annulus $A = \alpha \times [0,1]$ consists of 
$\bdry$-parallel arcs and possibly closed curves. }
\s

\begin{figure}
	{\epsfysize=2.5in\centerline{\epsfbox{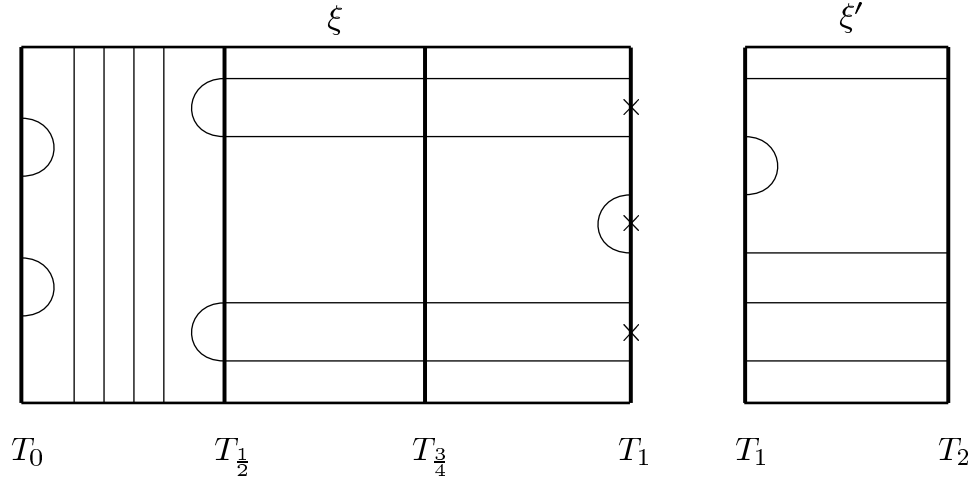}}}
	\caption{Templates.  The top and bottom are identified.}
	\label{newtemplate}
\end{figure}

Let us mark the 
components of $\Gamma_{T_1}$ that are straddled by the $\bdry$-parallel 
arcs on $A$ (by $\times$ 
in Figure \ref{newtemplate}). 
Note that adding a bypass to $\xi$ across any marked 
curve will produce an overtwisted
structure, since the dividing set on the extended annulus will contain a 
closed trivial dividing curve.
Any bypass we discuss here is along an arc
connecting three different dividing curves, 
and we call it a bypass ``across'' the middle curve.

\begin{prop} \label{attach=dig} Let  $(T^2\times [0,1],\xi)$ be a rotative universally 
tight contact structure which satisfies $\#\Gamma_{T_1}>>2$ and the $\bdry$-Parallel 
Assumption, and let $(T^2\times [1,2],\xi')$ be a non-rotative layer corresponding to attaching 
a bypass across an unmarked dividing curve on $T_1$. There exists a contact isotopy
$\phi_t:(T^2\times [0,2],\xi\cup \xi')\hookrightarrow (T^2\times [0,2],\xi\cup\xi')$, 
$t\in[0,1]$,  rel $T_0$, where $\phi_0=id$ and $\phi_1(T^2\times[0,2])\subset T^2\times[0,1]$.
\end{prop}

\begin{figure}
	{\epsfysize=2in\centerline{\epsfbox{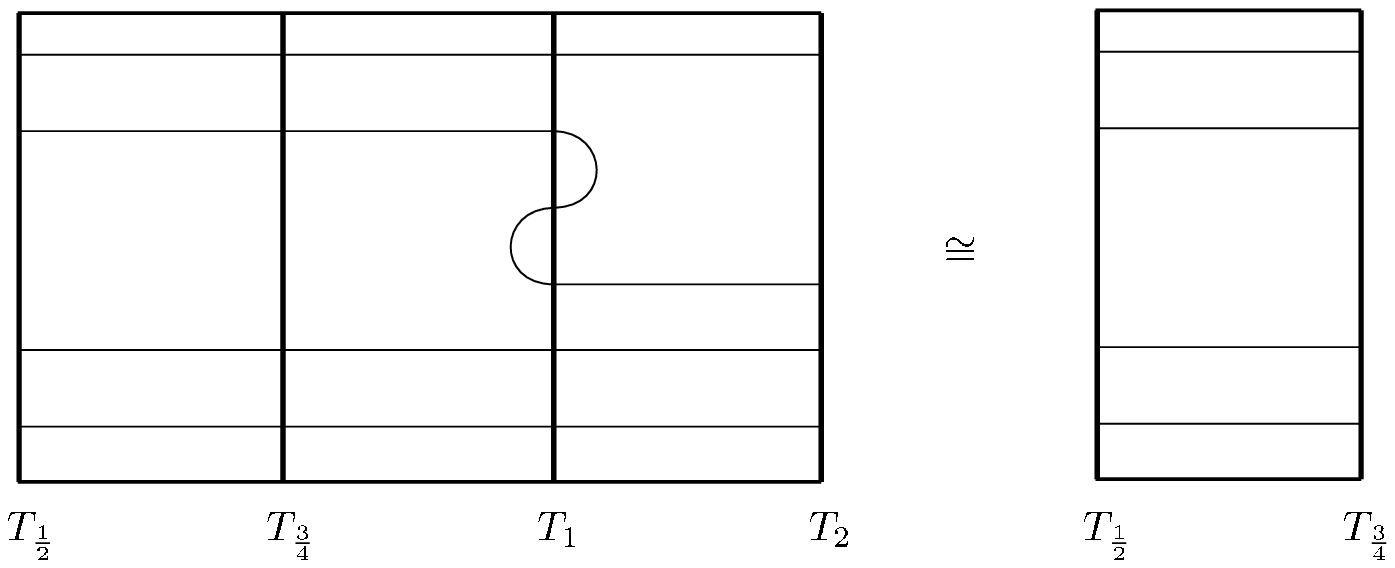}}}
	\caption{}
	\label{newattach=dig}
\end{figure}

\begin{proof} Observe that $(T^2\times [0,1],\xi)$ can be factored (after possibly isotoping 
relative to the boundary) 
into $(T^2\times[0,{1\over 2}], \xi) \cup (T^2\times[{1\over 2},1], \xi)$, 
where $T_{1/2}$ is  
convex, $\xi$ is rotative on $T^2\times [0,{1\over 2}]$,  non-rotative on 
$T^2\times[{1\over 2},1]$ and the dividing set induced on the annulus 
$A'' = \alpha \times [{1\over 2},1]$  consists of one $\bdry$-parallel arc 
and nonseparating arcs (see Figure \ref{newtemplate}). The factorization can be chosen so that 
the $\bdry$-parallel arc straddles one of 
the marked dividing curves adjacent to the unmarked curve across which the bypass  
of the layer $(T^2\times [1,2],\xi')$ is added.  
Note that $\#\Gamma_{T_{1/2}} = \#\Gamma_{T_1}  - 2$. 
Finally, consider $(T^2\times[{1\over 2},2],
\xi\cup \xi')$. It is isotopic to a product neighborhood
 $\xi|_{T^2\times[{1\over 2},{3\over 4}]}$ of $T_{1/2}$ through an isotopy which fixes  
$T_{1/2}$.  This is because a non-rotative contact structure is determined by the dividing set 
on a horizontal annulus (see Figure \ref{newattach=dig}).
\end{proof}

The same proof shows that digging out a single bypass along $A$ (removing a portion
like $T \times [{1\over 2},1]$ above) is equivalent to adding a bypass across an 
adjacent unmarked dividing curve.    
Therefore, we find that ``attaching a $T^2\times I$ layer'' is the same as
``digging out a $T^2\times I$ layer''.  

The Attach = Dig Principle applies more generally to a compact 3-manifold
$M$ with torus boundary $T$  assuming $M=M'\cup_{\bdry M'=T\times\{0\}} (T\times [0,1]) $, 
and $\xi$ is rotative on $T\times [0,1]$ with 
$\#\Gamma_{T\times \{1\}}>>2$ and satisfies the $\bdry$-Parallel Assumption. As in the torus 
case, we can use a horizontal annulus $A$ for $T\times [0,1]$ to mark dividing curves which are  
straddled by $\bdry$-parallel arcs in $\Gamma_A$. It follows from the previous proposition that 
attaching  a bypass across an unmarked dividing curve is equivalent to a dig. It takes a bit 
more work to argue that digging out a bypass is equivalent to attaching one.  Consider the 
factorization $M=M_B \cup (T\times[0,1])$, where $T \times [0,1]$ corresponds to
a  bypass attached to $T_1 \subset \bdry M$ on the inside along an arc connecting 
three different dividing curves (a ``short'' bypass).

First of all, we know that any ``short'' bypass on the inside must be 
happening across a marked curve. 
Indeed if there was an inside bypass across an unmarked curve 
there would be an obvious overtwisted disc
in the contact structure obtained by attaching the bypass on the outside 
across the same dividing curve. This would be  a contradiction since
we know that attaching such a bypass is equivalent to a 
dig and hence produces a tight structure 
(a subset of a tight structure is clearly tight).
We can now look at
the factorization $M=M_B \cup (T\times [0,1])$ and 
conclude that $(M_B,\xi)$ is isotopic to 
$(M_B \cup (T\times[0,1]) \cup (T\times [1,2]), \xi \cup \xi')$ 
with $(T\times [1,2],\xi')$ 
an outside bypass across an adjacent unmarked curve. 

The advantage of digging over attaching is manifest when we try to prove 
universal tightness by using the state 
transition method. We will also use the fact that 
digging can be viewed as attaching in proving a
gluing theorem along incompressible tori.

\subsection{Folding}       \label{section:folding}

\begin{figure}
	{\epsfysize=1in\centerline{\epsfbox{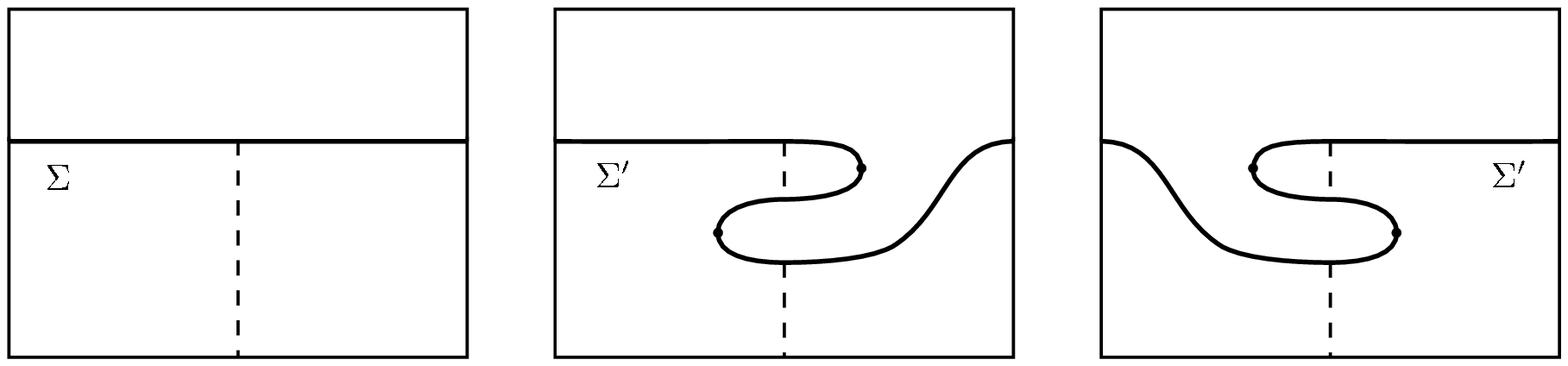}}}
	\caption{Folding.}
	\label{fold}
\end{figure}

We will now explain the process of {\it folding}, which plays a
crucial role in our proof.  A similar discussion can be found in \cite{H1}.
Let $\Sigma$ be a convex surface and
$\gamma$ be a {\it nonisolating} closed curve with $\gamma\cap
\Gamma_\Sigma=\emptyset$.
A closed curve $\gamma$ is called {\it nonisolating} if every component of
$\Sigma\backslash\gamma$  intersects $\Gamma_\Sigma$.  Now, the Legendrian realization 
principle (see \cite{H1}) states that if $X$ is a contact vector field transverse to $\Sigma$ 
and $\gamma$ is a nonisolating curve, then there is a $C^0$-small isotopy $\phi_t$, 
$t\in[0,1]$, supported in a neighborhood of $\Sigma$, for which 
\be
\item $\phi_0=id$,
\item $\forall t$, $\phi_t(\Sigma)\pitchfork X$,
\item $\forall t$, $\phi_t=id$ on $\Gamma_\Sigma$ and $\Gamma_{\phi_t(\Sigma)}=\Gamma_\Sigma$,
\item $\phi_1(\gamma)$ is Legendrian.
\ee
Therefore, by using the Legendrian realization principle and possibly modifying 
$\Sigma$, we may take $\gamma$
to be  Legendrian curve.  A slight strengthening of the Legendrian realization principle allows 
us to take $\gamma$ to be a {\it Legendrian divide}.  This means that all the points 
of $\gamma$ are tangencies, and there exists a local model $N=S^1\times
[-\varepsilon,\varepsilon]
\times [-1,1]$ with coordinates $(\theta,y,z)$ and 1-form
$\alpha=dz-yd\theta$ such that
$\Sigma\cap N=S^1\times [-\varepsilon,\varepsilon]\times \{0\}$ and
$\gamma=S^1\times \{0\}\times \{0\}$.        A {\it fold} is a
modification $\Sigma\sc \Sigma'$
where $\Sigma$ and $\Sigma'$ are isotopic and $\Sigma=\Sigma'$
outside $N$. Inside
$N$, we bend $\Sigma$ into $S^1$ times an S-shape so that $\Gamma_\Sigma'$
consists of $\Gamma_\Sigma$ outside $N$, together with two new
parallel curves isotopic to
$\gamma$.  The picture on the left-hand side of Figure \ref{fold}
represents $\Sigma$.
Here the rectangle represents
$[-\varepsilon,\varepsilon]\times[-1,1]$ (the horizontal direction
is the $y$-direction and the vertical direction is the
$z$-direction).  There are two choices for
folding (middle and right pictures), and they are not contact
isotopic (see next paragraph).
The dots (times $S^1$) represent the new dividing curves on $\Sigma'$.
Note that if we let $N'$ be the
``upper'' solid torus split by $\Sigma'$, the dividing curves of the
meridional disks (after Legendrian
realization to make the boundary of the meridional disks Legendrian) will
be given by the dotted lines
in Figure \ref{fold}.    Observe that we folded $\Sigma$ to obtain
$\Sigma'$ inside an $I$-invariant
neighborhood of $\Sigma$.  Hence, folding to increase the number of
dividing curves (in pairs)
is an operation that can be done to change $\Gamma_\Sigma$ {\it for
free} (at least for nonisolating curves).  In general, there are no
other operations (such as decreasing the number of dividing curves)
which can be performed on the
dividing sets without prior knowledge of the ambient manifold or at
least more global information.

\begin{figure}
	{\epsfysize=1in\centerline{\epsfbox{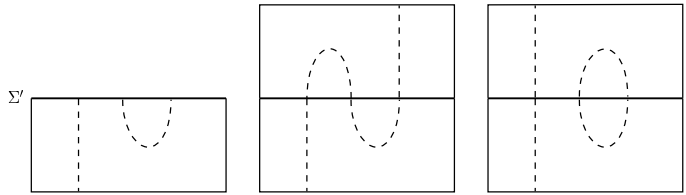}}}
	\caption{Templates.}
	\label{template}
\end{figure}

To distinguish the tight contact structures obtained by folding, we
apply the {\it template method}, which is similar to the method used 
in Section \ref{a=d} to identify the location of the $\bdry$-parallel arcs.  
Namely, we glue on various $S^1\times D^2$ with $S^1$-invariant
universally tight structures,
and distinguish the folds according to which templates
make the glued-up contact structure overtwisted.  For the single fold
above, we attach templates
depicted in Figure \ref{template}.   The top rectangles of the middle
and right diagrams are
the meridians of the templates.
Note that the middle attachment preserves tightness, but the right
attachment gives an overtwisted
contact structure.  On the other hand, if we took the mirror image of
the left-hand diagram (across a
vertical line), then the middle attachment will be overtwisted and
the right attachment will be universally
tight.  In a similar manner, all the folds are distinguishable by
attaching templates {\it from the outside}.

\section{Torus boundary case}      \label{section:boundary}
In this section we will first construct a suitable  tight contact structure
$\xi$ on $M$ by using Theorem \ref {gluing}.
Consider the situation where $\bdry M\not=\emptyset$
and $\bdry M=\cup_{i=1}^{l} T_i$, where $T_i$ are incompressible tori.
We choose a suitable surface $S$ as the first splitting surface in a well-groomed
sutured decomposition.

\begin{claim}
There exists a properly embedded, not necessarily connected,
well-groomed surface
$S$ which is (Thurston) norm-minimizing in
$H_2(M,\bdry M)$ and nontrivially intersects each $T_i$, that is, the image
under the boundary map
is a nonzero class in $H_1(T_i)$.
\end{claim}

\begin{proof}
For each $i=1,\cdots,l$, take $M(i)$ to be $M$ filled with solid tori
$S^1\times D^2$ along
each $T_j$ except for $T_i$.  This means $\bdry M(i)=T_i$.
Using the relative homology sequence and Lefschetz duality for the
manifold $M(i)$ and
its boundary $\bdry M(i)$, we find a class
$[S(i)]\in H_2(M(i),\bdry M(i))$ which is nonzero under the boundary map
$H_2(M(i),\bdry M(i))\stackrel{\bdry}{\rightarrow} H_1(M(i))$.
Represent it by a surface
$S(i)$ in $M(i)$.   Back on $M\subset M(i)$ let $S_i=S(i)\cap M$ and  consider
$[S(i)\cap T_j]$ for $j\not=i$ (assume
transversality).  If $[S_i\cap T_j]=0\in H_1(T_j)$ we can pair off
the intersections and get $S_i\cap T_j
=\emptyset$.    Now, take a suitable linear combination of the
$[S_i]$ to obtain $[S]\in
H_2(M,\bdry M)$ for which
$\bdry[S]$ is nonzero when restricted to each $H_1(T_i)$.  Finally
pick a norm-minimizing
representative of $[S]$.   The well-grooming is simply asking that
$T_i\cap \bdry S$ be
parallel curves oriented in the same direction. (If not already
well-groomed, we can pair off
oppositely oriented parallel curves without altering the norm.)
\end{proof}

To proceed with the construction, let $\Gamma_{\bdry M}$ consist of
pairs of parallel essential curves for each $T_i$, chosen so
that $\Gamma_{T_i}$
has nonzero intersection with $T_i\cap \bdry S$ on $T_i$ and that 
every component of
  $T_i\cap \bdry S$ is efficient with respect to $\Gamma_{\bdry M}$.
$\Gamma_{\bdry M}$ is taut, since each component of $R_\pm$ is an annulus.
The order of choosing things here might seem reverse from the way we 
study contact structures by
decomposing them along convex surfaces, and it is --- we are choosing 
the norm-minimizing surface
$S$ first,  and then choosing
$\Gamma_{\bdry M}$.  We are free to do that since we are constructing 
a contact structure by gluing,
rather than analyzing one that is already given to us. Let us
assume first that all the
$T_i$ in the sutured decomposition are torus sutures. We have chosen 
the next cutting surface $S$
to be norm-minimizing and well-groomed,  implying that the resulting 
sutured manifold
$M\backslash S$ is taut.
Now, choosing $\Gamma_{\bdry M}$ to have nonzero
geometric intersection with each $T_i$, together with the well-grooming of $S$,
guarantees that after splitting and rounding, the sutures (dividing curves) on
$\Gamma_{\bdry (M\backslash S)}$ are the same as would have been 
obtained from the sutured
decomposition in the case when all  the $T_i$ are torus sutures.
Theorem \ref{Gabai:hierarchy} ensures
a sutured manifold hierarchy where the first splitting is along $S$.
Note that Gabai's theorem also
ensures that each subsequent cutting is done along surfaces which are
well-groomed and have boundary. Therefore, by Theorem \ref{glueback}, we can construct
a universally tight contact structure $\xi$ on $M$ for which $S$ is convex with
$\Gamma_S$ $\bdry$-parallel.

Assume that we have fixed $\Gamma_{\bdry M}$ as above and constructed a specific universally 
tight $(M,\xi)$.
We will proceed to describe a family of contact structures
$\xi'_k$, $k\in \Z^+$, obtained by attaching $(T^2\times I, \xi_k)$ layers onto $M$ along 
$T=T_i$, where $\xi_k$ is defined as in Section \ref{section2-1}. 
Write $\xi_k^+=\xi_k$, $k\in \Z^+$, and let $\xi_k^-$ be given by $-\alpha_k$ 
(the 2-plane fields $\xi_k^+$ and $\xi_k^-$ are identical, but are oriented oppositely).
Consider the glued-up contact manifold $(M'=M\cup_T (T^2\times
I),\xi'_k)$, where we use some
element of $ SL(2,\Z)$ to glue $T^2\times \{0\}$ to $T$ so that the
dividing curves match up.

\begin{prop}    \label{step1}
$\xi'_k$ is
universally tight for precisely one of  $\xi_k^+$ or $\xi_k^-$.
\end{prop}

\proof Cut $M'$ open along $S'=S\cup ((\bdry S\cap T)\times I)$, where 
$T=T\times\{0\}$.
Denote the components of  $\bdry S$  by $\gamma_1,\cdots,\gamma_s$,
$\gamma_{s+1},\cdots,\gamma_t$, where $\gamma_i$, $i=1,\cdots,s$, are
parallel curves on
$T$ and the other $\gamma_i$ lie on different $T_j\not=T$.  Then
$S'\backslash S$ consists of annuli $\gamma_i\times I$, $i=1,\cdots,t$.
Let $A_1,\cdots,A_s$ be the annular components of
$(T\times\{1\})\backslash \bdry S'$ and
$A_{s+1},\cdots, A_t$ be the annular components of 
$\cup_{T_j\not=T} (T_j\times\{1\})\backslash \bdry S'$.  Recall that $S$ was constructed so 
that $T_j\cap S\not=\emptyset$ for all $T_j$.   By possibly 
reordering, we may assume that $\bdry A_i =(\gamma_{i-1}\times\{1\})\cup(\gamma_i\times\{1\})$, 
$i=1,\cdots,s$, where $\gamma_0\stackrel{def}{=} \gamma_{s}$.

First observe that Proposition \ref{horizontal} implies that, for a correct choice of
$\xi_k^+$ or $\xi_k^-$ (but not both, since in one case the 
$\bdry$-parallel components
will match to form a closed curve parallel to the boundary, and in 
the other to form
homotopically trivial circles, giving overtwisted disks) 
$\Gamma_{S'}$ will consist of
$\bdry$-parallel curves identical to those of
$\Gamma_S$,  together with
$k$ extra curves parallel to $\gamma_i$ for each $i=1,\cdots,s$.
$\Gamma_{\bdry(M\backslash S)}$ is isotopic to the collection of
curves which consists of a single core
curve from each $A_i$, $i=1,\cdots,t$; $\Gamma_{\bdry (M'\backslash
S')}$ is $\Gamma_{\bdry
(M\backslash S)}$ with $2k$ extra curves parallel to the core curve
of $A_i$, for each
$i=1,\cdots,s$.  

\begin{claim}
If every connected component of $S'$ has at least two boundary components, then
the contact structure $\xi'_k$ on $M'\backslash S'$ is the tight contact
structure $\xi$ on $M\backslash S$ with folds introduced.
\end{claim}

\begin{proof}[Proof of claim]
Observe that $M'\backslash S'=(M\backslash S)\cup((T \backslash (\cup_{i=1}^s 
\gamma_i))\times I)$. Let $\gamma'_i$, $i=1,\cdots, s$ be parallel copies of $\gamma_i$
on $S$ which do not intersect $\Gamma_S$, and let $B_i\subset S'$ be the annulus bounded by
$\gamma'_i$ and $\gamma_i\times\{1\}$.
Then take the annulus $\Sigma'\subset \bdry(M'\backslash S')$  to be
a union of $A_i$ and copies of $B_{i-1}$ and $B_i$.    Also let
$\Sigma\subset \bdry(M\backslash S)$ be an annulus parallel to $\Sigma'$ with
$\gamma'_i$, $\gamma'_{i+1}$ as boundary. 
We may use the
Legendrian realization principle to realize copies of 
$\gamma_i'$ and $\gamma_{i+1}'$ on $\Sigma$ (and $\Sigma'$) as Legendrian divides. 
It is at this point that the
assumption that every connected component of $S'$ has at 
least two boundary components is needed to conclude that the
nonisolating condition in the Legendrian realization principle is met.
Now, after rounding,  $\Sigma$
and $\Sigma'$ are disjoint  except along their common boundary $\gamma_i'\cup \gamma_{i+1}'$.  
Let $N$ be the solid
torus region bounded by
$\Sigma$ and $\Sigma'$.  The dividing curves on the meridian of $N$
will be as in Figure
\ref{meridian}.
\begin{figure}
	{\epsfysize=1.5in\centerline{\epsfbox{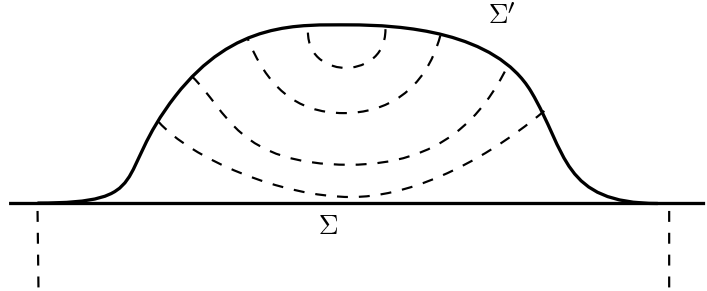}}}
	\caption{Meridian.}
	\label{meridian}
\end{figure}
With this explicit description, it is immediate that the attachment of $N$ is equivalent to a  
folding operation.  Finally, repeat this procedure with all annuli $A_i$, $i=1,\cdots, s$, in 
succession. 
\end{proof}

The assumption of the previous claim on the number 
of boundary components is a 
technical condition. It is satisfied by passing to a
large finite cover $\widetilde
{M'}$ which unwinds $S'$ to give $\widetilde {S'}$, which would have more than one boundary 
component even if $S'$ had one boundary component. The claim implies that there exists a 
contact embedding: $$(\widetilde {M'}\backslash \widetilde {S'},\widetilde
{\xi'_k})\hookrightarrow
(\widetilde M\backslash \widetilde S,\widetilde \xi).$$
Therefore, $\xi'_k$ on $M'\backslash S'$ is universally tight.

\begin{figure}
	{\epsfysize=2in\centerline{\epsfbox{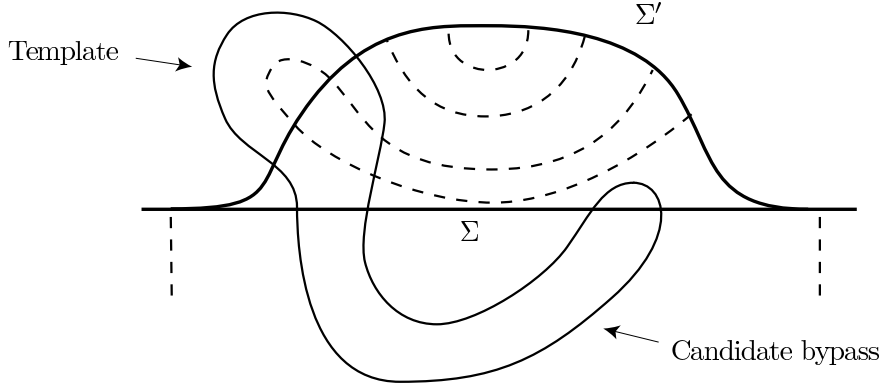}}}
	\caption{Attaching a template to obtain an overtwisted disk.}
	\label{templateattachment}
\end{figure}

Next, we need to deal with the problem of gluing $M'$ back along
$S'$.  The argument is almost identical
to the proof of the Gluing Theorem (Theorem \ref{gluing}) in the $\bdry$-parallel case.
Pass to a large finite cover $\widetilde{M}$ of
$M$, so that we have a large cover $\widetilde{S'}$ of
$S'$, and lift the candidate overtwisted disk $D$.  Here, every
bypass attachment will be either trivial or will
decrease the number of dividing curves.  A bypass along
$\widetilde{S'}$ which decreases the number of dividing
curves cannot exist by attaching a template.  See Figure
\ref{templateattachment}.
Therefore, $(M',\xi'_k)$ is universally tight.
\qed

\s\n
{\it Remark.}  The above proposition could have been proved using
Colin's Gluing Theorem
(Theorem \ref{pre-Lagrangian}) below.  Our point was to demonstrate a 
different perspective.

\s\n
{\it Change of notation.}  To avoid excessive use of primes, from now on we will refer to the 
contact manifold $(M',\xi_k')$ and the splitting surface $S'$ as $(M,\zeta_k)$ and $S$, 
respectively.   Also write $\bdry M=\cup_{i=1}^l T_i$

\begin{prop}    \label{prop:distinct}
The universally tight contact structures $\zeta_k$ are nonisotopic rel boundary.
\end{prop}

\begin{proof}
Assume  without loss of generality that $S$ has collared
Legendrian boundary. Let $T$ be the connected component of $\bdry M$ in whose direction twisting 
was added.  Now let $\delta:[0,1]\rightarrow S$ be a properly embedded arc on $S$ satisfying 
the following:

\be
\item[(a)] $\delta$ is essential (i.e., not $\bdry$-compressible).
\item[(b)] The endpoints $\delta(0)$, $\delta(1)$ lie on the {\it Legendrian divides} of the 
collared Legendrian boundary of $S$.   
\item[(c)] At least one of $\delta(0)$, $\delta(1)$ (both if 
possible) lies on $T$.   
\item[(d)] $\delta$ has minimal geometric intersection
$|\delta\cap \Gamma_{S}|$ amongst arcs on $S$ isotopic to $\delta$ with
endpoints on the same component of $\bdry S$ as $\delta$. 
\item[(e)] $\delta$ is Legendrian. (This is made possible by the Legendrian realization 
principle.) 
\ee
                                  
Since $\Gamma_S$ consists of 
$\bdry$-parallel arcs and closed curves parallel to $\bdry S$, any choice of
$\delta$ will have no intersections with the $\bdry$-parallel dividing arcs and minimally 
intersect the closed dividing curves parallel to the boundary, i.e., $k$ times if 
only one of $\delta(0)$, $\delta(1)$ lie on $T$,  
and $2k$ times if they both do. Define $\mathcal{L}_\delta$ to be the 
set of Legendrian arcs $d$ satisfying the following:
\be
\item There exists a properly embedded convex surface $S'\supset d$ with 
collared Legendrian boundary. 
\item $\bdry S'$ is efficient with respect to $\Gamma_{\bdry M}$.
\item The endpoints of $d$ lie on Legendrian divides of the collared Legendrian boundary 
of $S'$. 
\item There exists an isotopy (not necessarily a contact isotopy) which sends 
$S$ to $S'$ and $\delta$ to $d$, where $\bdry S$ is isotoped to $\bdry S'$ along $\bdry M$. 
\ee

We are now ready to define our invariant which distinguishes the $\zeta_k$.
For any $d\in \mathcal{L}_\delta$ we define the {\it twisting number}
$t(d)$ to be
$-{1\over 2}|d\cap \Gamma_{S'}|$, where $S'$ is the convex surface
containing $d$ in the definition of $\mathcal{L}_\delta$ above.
Define the {\it maximum twisting} for $\mathcal{L}_\delta$ to be
$$t(\mathcal{L}_\delta)=\max_{d\in \mathcal{L}_\delta} t(d).$$

We claim that $t(\mathcal{L}_\delta)=t(\delta)$.
Assume on the contrary that there exists a Legendrian arc $d\in \mathcal{L}_\delta$
with $t(d)=t(\delta) + 1$, sitting on the convex surface $S'$. There exists an isotopy 
$\phi_t:[0,1]\rightarrow M$, $t\in[0,1]$, such 
that $\phi_0=\delta$ and $\phi_1=d$. 
Since  $\bdry S$ and $\bdry S'$ are efficient with respect to $\Gamma_{\bdry M}$ and are 
isotopic on $\bdry M$ through efficient curves, it is possible to ``slide'' $S$ along $\bdry M$ 
without altering the isotopy class of $\Gamma_S$, and we may assume that $\bdry S=\bdry S'$ and 
$\phi_t(0)$, $\phi_t(1)$ are fixed. As usual, pass to a large {\it finite} cover $\widetilde M$ 
of $M$ to unwrap $S$ to $\widetilde{S}$.

\begin{claim}
For a large enough finite cover $\widetilde {M}$, there exists a
lift/extension of $\phi_t:[0,1]\rightarrow
M$, $t\in[0,1]$, to $\Phi_t:\widetilde{S}\rightarrow \widetilde {M}$,
where $\Phi_0$ is the identity and the support of the
isotopy $\Phi_t$ is compact and is contained inside an embedded 3-ball.
\end{claim}

\begin{proof}[Proof of claim]
The key property we use here is that the isotopy $\phi_t$,
$t\in[0,1]$, can be lifted to a large finite cover
$\widetilde {M}$ so that the trace of the lift
$\widetilde{\phi}_t:[0,1]\rightarrow \widetilde {M}$
is contained in a 3-ball $B^3$.  This follows from the fact that
$M$ is Haken, which implies
that  $M$ has universal
cover $\R^3$   and $\pi_1(M)$ is residually finite.
Let $\widetilde X_t$, $t\in[0,1]$, be a time-dependent vector field
along $\widetilde{\phi}_t([0,1])$
which induces the isotopy $\widetilde{\phi}_t$.
Damp $\widetilde X_t$ out away from $\widetilde{\phi}_t([0,1])$ and
extend it to all of $\widetilde {M}$.  Then $\cup_{t\in[0,1]}
supp (\widetilde X_t)\subset B^3$, where the $supp(X)$ is the support
of a vector field (i.e., $p\in M$
such that $X(p)\not=0$).  Let $\Phi_t$, $t\in[0,1]$, be the flow
corresponding to $\widetilde X_t$.
It is supported inside $B^3$ from our construction.
\end{proof}

By assumption, $|t(d)|$ is strictly smaller
than $k$ (or $2k$). If we can prove that, for each component $C$ of $\bdry S'$ which is on $T$, 
there exist at least $k$ closed curves in $\Gamma_{S'}$ parallel to $C$, we have a 
contradiction. In fact it is clearly enough to obtain the contradiction for any large enough
cover.
Take the cover $\widetilde {M}$ of $M$ given by the previous claim, and
use the technique of {\it isotopy discretization},
which first appeared in Colin \cite{Co97}, and is a method which
works extremely well in the context of convex surfaces and state
transitions (due to bypasses).  We split
$[0,1]$ into small time intervals $t_0=0<t_1<\cdots< t_{n-1}<t_n=1$,
$n>>0$, so that all the  $\Phi_t$ are
are all ``sufficiently close'' if $t\in [t_i,t_{i+1}]$, in the following sense:
there exists an embedding $\psi_i:\widetilde{S}\times I\rightarrow
\widetilde {M}$ such that
$\Phi_t(\widetilde{S})\subset \psi_i(\widetilde{S}\times I)$ and
$\pi_1\circ \psi_i^{-1}\circ \Phi_t(\widetilde{S})$ is a submersion
for all $t\in[t_i,t_{i+1}]$.
(Here $\pi_1:\widetilde{S}\times I\rightarrow \widetilde{S}$ is the
first projection.)
This implies (by Giroux's convex movie approach in \cite{Gi99b} and the equivalence 
of Giroux's `retrogradient switch' with the notion of a bypass) that to get from 
$\Phi_0(\widetilde{S})=\widetilde{S}$ to $\Phi_1(\widetilde{S})$ we move via  a sequence of 
bypass moves. We may need to
take an even larger cover of $\widetilde {M}$ to expand $\widetilde
{S}$ (without changing names)
and extend $B^3$ so that we can ``localize'' all
the bypasses, i.e.,  all the allowable bypasses are either trivial or
reduce the number of dividing curves.
The bypasses which potentially give us trouble are the bypasses which
increase the number of dividing curves.
Such a bypass occurs only when the bypass ``wraps around'' a closed curve which
is homotopically nontrivial ---
by unwrapping in this direction in a large enough cover, we may avoid 
a dividing curve increase.
More precisely, since $S$ is incompressible in $M$, $\pi_1(S)$ injects into 
$\pi_1(M)$, and we may unwrap $S$ in the direction of any homotopically nontrivial
curve in a finite cover of
$M$ using residual finiteness.

We can show that there are no bypasses reducing the number of dividing curves
as in the proof of Proposition \ref{step1}. Therefore,
we keep attaching trivial
bypasses, and eventually find a contact isotopy of $\widetilde{S}$ with
$\Phi_1(\widetilde{S})$, which contradicts the assumption that
$t(d)> t(\delta)$. \end{proof}

More generally, let $\zeta_{k_1,\cdots,k_l}$ 
be the universally tight contact structure obtained by attaching appropriate $(T^2\times I, 
\xi_{k_i}^\pm)$ layers  ($i=1,\cdots, l$) to each component $T_i$ of $\bdry M$.    

\begin{prop} The universally tight contact structures $\zeta_{k_1,\cdots,k_l}$ 
are distinct up to isotopy rel boundary. 
\end{prop}

\section{Closed case}      \label{section:closed}

\subsection{Topological preliminaries}

\subsubsection{Residual finiteness} \label{section:residual} A 
group $G$ is said to be {\it residually finite} if   given any nontrivial element $g\in G$ 
there exists a finite index {\it normal} subgroup $H\subset G$ which does not contain $g$.   
When $G$ is the fundamental group $\pi_1(M)$ of a manifold $M$, the residual finiteness of $G$ 
is equivalent to the following:  If $K$ is a compact subset of the universal cover $\overline 
M$ of $M$, then there exists a finite cover $\widetilde M$ of $M$ for which the projection 
$\overline M\rightarrow \widetilde M$ is injective on $K$.  Note that since $\pi_1(\widetilde 
M)$ can be taken to be {\it normal} in $\pi_1(M)$, $\widetilde M$ is a {\it Galois} (or {\it 
regular}) cover. We list some facts about residual finiteness (which can be found in 
\cite{Bo00}):

\be
\item A subgroup of a residually finite group $G$ is residually finite.
\item (Hempel \cite{He}) If $M$ is Haken, then $\pi_1(M)$ is  residually finite.
\ee

Let $H$ be a subgroup of a group $G$.  We say $H$ is {\it separable} in $G$ if for any element 
$g\in G - H$ there exists a finite index subgroup $K$ in $G$ satisfying $H\subset K$ and 
$g\not\in K$.   The following is a result of Long and Niblo \cite{LN}:

\begin{thm}   \label{separability}
Let $M$ be a closed oriented Haken 3-manifold.  If  $i: T^2 \rightarrow M$ is an 
incompressible embedded torus, then  $i_*(\pi_1(T^2))$ is a separable subgroup in $\pi_1(M)$.
\end{thm}

\s\n
Let $M$ be a compact manifold with an incompressible torus $T$.  We explain what 
it means to {\it expand} $T$ and obtain a finite cover $\widetilde M$ of $M$.  By property (1) 
above, if  we choose a finite collection $\mathcal{C}$ of elements of $\pi_1(T)=\Z^2$, there 
exists a finite index normal subgroup $H$ which avoids $\mathcal{C}$.  If $\pi_1(T)$ has basis 
$\alpha,\beta$, and we take $\mathcal{C}=\{m\alpha+n\beta | 0<m^2+n^2\leq R\}$ for some large 
$R>0$, we obtain a finite Galois cover $\widetilde M$ of $M$ which {\it expands} $T$.

\subsubsection{Torus decompositions}  The following is the Torus Decomposition Theorem (also 
often called the JSJ decomposition), due to Jaco-Shalen \cite{JS} and Johannson \cite{Jo}.  
\begin{thm}    \label{JSJ} 
Let $M$ be a closed, oriented, irreducible 3-manifold.  Then, up to isotopy, there exists a 
unique collection $\mathcal{T}$ of disjoint incompressible tori $T_1,\cdots,T_m$ which 
satisfies the following: 
\be
\item Each component of $M\backslash \cup_{i=1}^m T_i$ is atoroidal (i.e., any 
incompressible torus can be isotoped into a boundary torus), or else is Seifert fibered. 

\item $\mathcal{T}$ is minimal in the sense that (1) fails when 
any $T_i$ is removed.
 \ee
\end{thm}

Most atoroidal components can be endowed with a 
hyperbolic structure with cusps, except for the following potential components which are 
Seifert fibered: (1) Seifert fibered spaces with $m$ singular fibers over a base $\Sigma$ which 
is $S^2$ with $n$ punctures, where $m+n\leq 3$.    (2) Seifert fibered space with $m$ singular 
fibers over a base $\Sigma$ which is $\R\P^2$ with $n$ punctures, where $m+n\leq 2$.

A clear account of the Torus Decomposition Theorem appears in \cite{Ha}.  The following facts, 
also found in \cite{Ha} are useful. 

\begin{prop}  \label{horizontal-vertical}
Let $M$ be a connected, compact, oriented, irreducible Seifert fibered space.  Then any 
incompressible, $\bdry$-incompressible properly embedded surface $S\subset M$ is isotopic to a 
surface which is either vertical (union of regular fibers), or horizontal (transverse to all 
fibers). \end{prop}

\begin{lemma}    \label{annulus}
Let $M$ be a connected, compact, oriented, irreducible, and atoroidal.  If $M$ contains an 
incompressible, $\bdry$-incompressible annulus meeting only torus components of $\bdry M$, then
$M$ is a Seifert fibered space. 
\end{lemma}

\subsubsection{Diffeomorphism groups of 3-manifolds}  We also need to recall the following fact 
regarding the mapping class group $\pi_0(\mbox{Diff}(M))$ of a Haken 3-manifold 
$M$.

\begin{thm}    \label{diffeo}   Let $M$ be an oriented, compact Haken 3-manifold and 
$\text{Diff}(M)$ its group of diffeomorphisms.  Then 
$\pi_0(\text{Diff}(M))$ is finite modulo Dehn twisting along possible incompressible tori and 
incompressible, $\bdry$-incompressible annuli. 
\end{thm}

For a 3-manifold, a {\it Dehn twist} is a diffeomorphism $\phi$ which is locally given on 
$T^2\times I=(\R^2/\Z^2)\times [0,1]$ or $[0,1]\times (\R/\Z)\times [0,1]$ with 
coordinates $(x,y,t)$  by: $\phi:(x,y,t)\mapsto (x,y+t,t)$.

\subsection{Gluing along incompressible tori}

The following result, due to Colin \cite{Co99a}, gives a sufficient condition for gluing tight 
contact manifolds along tori: 

\begin{thm}[Colin] \label{pre-Lagrangian}
Let $M$ be an oriented, irreducible 3-manifold and $T\subset M$ an 
incompressible torus.
Consider a contact structure $\xi$ for which $T$ is pre-Lagrangian 
(linearly foliated)  and
$\xi|_{M\backslash T}$ is universally tight.  Then $\xi$ is 
universally tight on $M$.
\end{thm}

From the perspective of convex decomposition theory, this is wholly 
unsatisfactory, since the
gluing surface is not convex!  We therefore formulate a more 
cumbersome variant of Colin's
Gluing Theorem which is phrased and proved entirely using convex 
surface theory.

\begin{thm}[Variant of Colin's Theorem] \label{variant}
Let $M$ be an oriented, irreducible 3-manifold and $T\subset M$ an 
incompressible torus.
Consider a contact structure $\xi$ for which $T$ is convex, $\xi$
is universally tight when restricted to $M\backslash T$, there exists 
a toric annulus layer
$N=T\times[-1,1]\subset M$, where $T=T\times\{0\}$, $\xi|_N$ is 
universally tight and
$\xi|_{T\times [-1,0]}$ and $\xi|_{T\times [0,1]}$ are both rotative. 
Then $\xi$ is
universally tight on $M$. \end{thm}

\begin{proof}
Recall that, according to the classification of tight contact 
structures on toric annuli (see
\cite{H1}), the universal tightness of $\xi|_N$ simply means that it 
is obtained by taking
the standard rotative contact structure $\sin(2\pi z)dx+\cos(2\pi 
z)dy=0$ on $T^2\times \R$
with coordinates $((x,y),z)$, truncating for a certain interval 
$[z_0,z_1]\subset \R$,
making the boundary convex with $2$ dividing curves each, and finally 
folding to attach a
non-rotative layer.   Without loss of generality, we may take both
$\xi|_{T\times [-1,0]}$ and $\xi|_{T\times [0,1]}$ to satisfy the $\bdry$-Parallel 
Assumption with $\bdry$-parallel arcs along $T_0$. 

We will prove that $(M,\xi)$ is tight by
contradiction.  The proof of tightness of finite covers of
$M$ is identical, since any finite cover of $(M,\xi)$ also satisfies 
conditions of Theorem
\ref{variant}.
This implies the universal tightness of $\xi$, since Haken
3-manifolds have residually finite fundamental group, so any 
overtwisted disk that would exist in
the universal cover could be 1-1 projected into a finite cover.

Assume $(M,\xi)$ is overtwisted with an overtwisted disk $D\subset M$.
As before, we look at a large finite cover
$\pi:\widetilde M\rightarrow M$ that sufficiently expands $T$.  Let
$\widetilde T_i$ be connected components of $\pi^{-1}(T)$ and $\widetilde D$ be a lift of $D$ 
to $\widetilde M$.  In $\widetilde M$ we may {\it extricate} $\widetilde D$ 
from $\widetilde T_i$ essentially without any penalty, as we shall see.
As in the proof of Theorem \ref{gluing}, we push $\widetilde D$ across 
$\widetilde T_i$ in stages.  We can eliminate components of
$\widetilde D\backslash (\cup_i
\widetilde T_i)$ by pushing across balls, and this is equivalent 
to a sequence of bypass moves.  If we expand
$T$  sufficiently, i.e., we take a large enough cover, the only bypass 
operations are either (1) trivial or (2) reduce
$\#\Gamma_{\widetilde {T_i}}$.   In particular, there are no dividing 
curve increases or
changes in slope.

Let $\widetilde T$ be a fixed $\widetilde T_i$.                                                                           
Denote by $\widetilde T\times [-1,1]$ the connected component of $\pi^{-1}(T\times[-1,1])$ 
where $\widetilde T=\widetilde T\times \{0\}$.  Suppose the 
bypass attached onto $\widetilde T$ from the $\widetilde T\times 
[0,1]$ side decreases $\#\Gamma_{\widetilde T}$.  Let $\widetilde
T'$ be the new torus we obtain after the bypass attachment.  Then $(\widetilde 
M\backslash \widetilde T',\pi^*(\xi))$ is tight by the Attach = Dig Principle  (Section 
\ref{a=d}). Repeated application of the Attach = Dig Principle proves 
Theorem \ref{variant}. 
\end{proof}

\subsection{Proof of Theorem \ref{infinite}} 

Let $M$ be a closed manifold and $T\subset M$ an oriented incompressible torus.
The cut-open manifold $M\backslash T$ will have two torus boundary components which we
denote by $T_+$ and $T_-$;  $T_+$ is the component where the boundary orientation agrees 
with the orientation induced from $T$, and $T_-$ is the component where they are opposite. 
$T$ can be either separating or nonseparating.
If $T$ is separating, then $M\backslash T$ will consist of two
components $M_1$, $M_2$, with
$\bdry M_1=T_+$, $\bdry M_2=-T_-$.
If $T$ is nonseparating, $M\backslash T$ will consist of one component
$M_1$ with $\bdry M_1=T_+ - T_-$.
In either case, initially choose $\Gamma_{\bdry (M\backslash T)}$ with 
$\Gamma_{T_+}=\Gamma_{T_-}$ so that 
splitting surfaces $S_i$, chosen below, intersect
 $T_+$ in a collection of parallel essential curves
which are not parallel to
$\Gamma_{T}$ (and likewise for $T_-$).

Construct universally tight contact structures $\zeta_{k_+,k_-}$ on $M\backslash T$ as in 
Section \ref{section:boundary} by (1) performing a convex decomposition corresponding to a 
sutured manifold decomposition and regluing as in the proof of 
Theorem \ref{thm1} and (2) adding a $\pi k_+$-twisting layer to $T_+$ and a $\pi k_-$-twisting 
layer to $T_-$ as in Proposition \ref{step1}.   For the correct parity of $k=k_++k_-$, Theorem 
\ref{variant} gives a universally tight contact structure on $M$ which we call $\zeta_k$.

\begin{thm}   \label{isotopy}
If $M$ is a closed, oriented 3-manifold with an incompressible torus,
then there exist infinitely many
universally tight contact structures up to isotopy.
\end{thm}

\begin{proof}
We define the invariant which will allow us to distinguish the $\zeta_k$ up to 
isotopy.  First we choose a suitable isotopy class of closed curves.  If $T$ separates, then 
take the next splitting surfaces
$S_i$ of $M_i$, $i=1,2$.  Since $S_i$ cannot be a disk with boundary
on $T$, we may take
arcs $\gamma_i\subset S_i$ with endpoints on $T$ which are not
$\bdry$-compressible
on $S_i$ as well as in $M_i$.  Glue $\gamma_1$ and $\gamma_2$ to
obtain $\gamma$.
If $T$ is nonseparating and there is an (incompressible) splitting
surface $S$ which spans between
$T_+$ and $T_-$, then take $\gamma_1$ to be some arc on $S$ with endpoints on
$T_+$ and $T_-$ and glue the endpoints to obtain $\gamma$.  If there
is no single surface
$S$, then take $S_1$ to have boundary component(s) on $T_+$ and $S_2$
to have boundary
components on $T_-$, find $\gamma_i$, $i=1,2$, as before, and glue to
get $\gamma$.
Let $\mathcal{L}_\gamma$ be the set of Legendrian curves which are
isotopic to $\gamma$
(but not necessarily Legendrian isotopic).

Assume without loss of generality that we have two surfaces $S_1$ and $S_2$.    Recall 
$\bdry S_i$, $i=1,2$, consist of efficient parallel oriented curves on $T$ which have nonzero 
geometric intersection with $\Gamma_T$.  We construct a surface $S_1+S_2$ by performing a band 
sum at each (geometric) intersection of $\bdry S_1$ and $\bdry S_2$.  Whether we do a positive 
band sum or a negative band sum depends on the following.  In order to define the band sum, 
temporarily assume $\Gamma_T$, $\bdry S_1$, $\bdry S_2$ are all linear on $T$, with slopes 
$\infty$, $s_1$, $s_2$, respectively. (That is, for the time being, forget the fact that $\bdry 
S_1$ and $\bdry S_2$ are Legendrian, and treat them just as curves on $T$.)  Thicken $T$ to 
$T\times [-\varepsilon,\varepsilon]$ so that $\bdry S_1\subset T\times\{\varepsilon\}$ and 
$\bdry S_2\subset T\times\{-\varepsilon\}$.  Each band intersects $T\times\{t\}$, 
$t\in[-\varepsilon,\varepsilon]$, in a linear arc of slope $s(t)$, which interpolates between 
$s_1$ and $s_2$ and is never $\infty$.      In the case $s_1=s_2$, take suitable multiples of 
$S_1$ and $S_2$ (still call them $S_1$ and $S_2$) so $\bdry S_1=\bdry S_2$, and call this 
$S_1+S_2$.

We now define the framing for $\gamma$ to be one which comes from $S_1+S_2$.  
Here we are assuming without loss of generality that $\gamma_i$, $i=1,2$, have endpoints on 
the same band.  A natural Legendrian representative of $\gamma$ will have subarcs $\gamma_i$, 
$i=1,2$,  where the endpoints of $\gamma_i$ lie on half-elliptic points on $\bdry S_i$, which 
are also on $\Gamma_T$.  This means that the actual Legendrian curves $\bdry S_1$ and 
$\bdry S_2$ intersect in a tangency at the endpoints of $\gamma_i$.  
As before, define the {\it maximal twisting} $t(\mathcal{L}_\gamma)$ to be the maximum 
twisting number of Legendrian curves in $\mathcal{L}_\gamma$. The following proposition proves 
Theorem \ref{isotopy}. \end{proof}

\begin{prop}     \label{twisting-number}
$t(\mathcal{L}_\gamma)=-{1\over 2}(k+1)$
or $-(k+1)$, depending on whether $\gamma$ intersects $T$ once or twice.
\end{prop}

\begin{proof}
It is clear that $t(\mathcal{L}_\gamma)\geq -{1\over 2}(k+1)$ or
$-(k+1)$, for these bounds are achieved by $\gamma$ on $S_1+S_2$.  We need to show that 
$t(\mathcal{L}_\gamma)\leq -{1\over 2}(k+1)$ or $-(k+1)$.
We will assume that there exist two surfaces $S_1$ and $S_2$ --- the
other case is similar.
Suppose $\gamma'$ is a Legendrian curve in $\mathcal{L}_\gamma$ with
$t(\gamma')=-(k+1)+1$.   There exists an isotopy $\phi_t:S^1\rightarrow M$, $t\in[0,1]$, for 
which $\phi_0(S^1)=\gamma$ and $\phi_1(S^1)=\gamma'$.
We pass to a large cover of $T$ to remove extra intersections of
$\gamma'$ with $T$.

In order to reduce from the case where $M$ is a closed manifold with 
incompressible torus  $T$
to the manifold-with-torus-boundary case, we construct an auxiliary 
contact manifold
$(W,\widehat\zeta)$ which has the following properties:

\be
\item $W$ is a large finite cover of $M$ (not necessarily Galois) and $\widehat\zeta$ the 
pullback of $\zeta_k$. 

\item $\phi_t:S^1\rightarrow M$, $t\in[0,1]$, lift to closed curves
$\Phi_t:S^1\rightarrow W$. In other words, $\Phi_t(S^1)$ maps 1-1 down to 
$\phi_t(S^1)$ under the covering projection $\rho:W\rightarrow M$.  

\item $T$ has been expanded sufficiently on $W$ so that the intersection of 
$\cup_{t\in[0,1]} \Phi_t(S^1)$ with any component $\widehat T_i$ of 
$\rho^{-1}(T)$ is (either empty or) contained  inside a ``small'' 
disk $D_i\subset \widehat T_i$ which does not intersect most of the 
dividing curves of $\widehat T_i$.

\item If we denote $\widehat\gamma=\Phi_0(S^1)$ and 
$\widehat\gamma'= \Phi_1(S^1)$, then $t(\gamma)=t(\widehat\gamma)$ and
$t(\gamma')=t(\widehat\gamma')$, where $t(\widehat\gamma)$ and $t(\widehat\gamma')$ are 
measured with respect to preimages $\widehat S_1$ and $\widehat S_2$ of $S_1$ and $S_2$.

\ee

The construction of $W$ will be done in the next section.  In the meantime we continue with the 
proof of Proposition \ref{twisting-number}, assuming $W$ has been constructed. 
$\Phi_t$ gives rise to
a time-dependent vector field $X_t$ along $\Phi_t(S^1)$; extend 
$X_t$ and hence $\Phi_t$ to all of $W$ by damping out outside a small neighborhood of 
$\Phi_t(S^1)$. Then the support $supp(X_t)$ of $X_t$ is contained in 
$N(\Phi_t(S^1))$.   
For each $\widehat T_i$, $\Phi_t(\widehat T_i)$ remains constant outside a disk $D_i\subset 
\widehat T_i$ which we may take to be convex with Legendrian boundary. 
If we discretize this isotopy into small time intervals $t_0=0< t_1<\cdots<t_n=1$, then 
each step $\Phi_{t_k}(\widehat T_i)$ to $\Phi_{t_{k+1}}(\widehat T_i)$ corresponds to
attaching a bypass which either is trivial or decreases the number of dividing curves while 
keeping the slope unchanged. Using the ``attach = dig'' principle, we find that the contact 
manifold $(W\backslash (\cup_i \Phi_1(\widehat T_i)),\widehat \zeta)$ is obtained from 
$(W\backslash (\cup_i \widehat T_i),\widehat \zeta)$ by attaching dividing-curve-decreasing 
bypasses or ``unfolding''.  Observe that $\widehat\gamma'$ intersects each $\Phi_1(\widehat 
T_i)$ at most once, say at $p_i$.  For each $p_i$, we may assume $p_i\in\Gamma_{\Phi_1(\widehat 
T_i)}$ and the isotopic copies $\widehat S_j'$ of $\widehat S_j$, $j=1,2$, which contain 
$\widehat \gamma'$ are convex with {\it efficient Legendrian boundary} (i.e., efficient with 
respect to the torus containing the boundary component).   In order to compare the twisting 
number of $\widehat\gamma'$ cut open along $\cup _i\Phi_1(\widehat T_i)$    and that of 
$\widehat \gamma$ cut open along $\cup_i \widehat T_i$, we need to modify the cut-open 
$\widehat\gamma$ slightly by attaching unfolding layers onto $W\backslash (\cup_i\widehat T_i)$ 
and extending this cut-open arc slightly on the attached toric annulus without modifying the 
twisting number.  Once we do this, we are comparing (the modified) $\widehat\gamma$  with 
$\widehat \gamma'$, both on $W\backslash (\cup_i \Phi_1(\widehat T_i))$, which reduces the 
problem to the torus boundary case. \end{proof}

\begin{proof}[Completion of proof of Theorem \ref{infinite}.]
Theorem \ref{isotopy} and its proof imply Theorem \ref{infinite}, once we look at the mapping 
class group $\pi_0(\mbox{Diff}(M))$ of the closed toroidal manifold $M$. 
Our argument is a little different from that of Colin \cite{Co99b} in that we do not bound the 
{\it torsion}, and instead use facts about the mapping class group to pass from ``infinitely 
many isotopy classes'' to ``infinitely  many isomorphism classes''.     Let $T_1,\cdots,T_m$ be 
the incompressible tori as in the Torus Decomposition Theorem and $M_1,\cdots,M_n$ the 
connected components of $M\backslash (\cup_{i=1}^m T_i)$.  Apply Gabai's sutured manifold 
decomposition to each $M_j$, obtain the corresponding universally tight contact structure on 
$M_j$, and add appropriate $(T^2\times I, \xi_{k_i})$ along each $T_i$.  This is 
$\zeta=\zeta(k_1,\cdots,k_m)$.  

Note that any element of $\mbox{Diff}(M)$ must fix the isotopy class of $\cup_{i=1}^m T_i$, but 
may permute the $T_1,\cdots,T_m$ (as well as $M_1,\cdots,M_n$).  Since there are only finitely 
many such permutations, we may assume that the elements of $\mbox{Diff}(M)$ we consider fix
each connected component $M_j$ (as well as $T_i$).
The hyperbolic components $M_j$ present no problem, since $\pi_0(\mbox{Diff}(M_j))$ is finite 
modulo Dehn twists on the boundary, according to Theorem \ref{diffeo}.   This means that, for 
any compact surface $S$ with boundary on $\bdry M_j$, there are finitely many isotopy classes 
of images of $S$ under $\mbox{Diff}(M_j)$, up to isotoping $\bdry S$ along $\bdry M_j$.  Let 
$M_j$ and $M_{j+1}$ be adjacent hyperbolic components with common boundary $T_i$, and $S_j$, 
$S_{j+1}$ be norm-minimizing surfaces in $M_j$, $M_{j+1}$ with boundary on $T_i$.   Both $S_j$ 
and $S_{j+1}$ are not disks or annuli (by Lemma \ref{annulus}).  Then take curves $\gamma_j$ on 
$S_j$ and $\gamma_{j+1}$ on $S_{j+1}$ to form a closed curve $\gamma$ as in the proof of 
Theorem \ref{isotopy}.  Since $S_j$, $S_{j+1}$ are not disks or annuli, we may take $\gamma_j$, 
$\gamma_{j+1}$ to be nontrivial curves which begin and end on the same boundary component of 
$S_j$, $S_{j+1}$, respectively.   Now consider $t(\mathcal{L}_\gamma)$ as well as 
$t(\mathcal{L}_{\phi(\gamma)})$, for all $\phi\in \pi_0(\mbox{Diff}(M))$.     Since there are 
finitely many isotopy classes $\phi(\gamma)$ modulo Dehn twists along tori, and Dehn twisting 
does not change $t(\mathcal{L}_\gamma)$, there exists a finite number 
$-k_\gamma=\inf_{[\phi]\in \pi_0(\mbox{\tiny{Diff}}(M))} t(\mathcal{L}_{\phi(\gamma)})$.  In a 
manner similar to Proposition \ref{twisting-number}, we can prove that 
$t(\mathcal{L}_\gamma)=-(k_i+1)$.   (Note that our current situation is not quite the same as 
the situation in Proposition \ref{twisting-number}, but the proof goes through, if we first 
extricate an isotopic copy of $\gamma$ from components $M_s$ with $s\not= j, j+1$.  This 
extrication can be done as before without penalty, by passing to a large finite cover.)  If we 
inductively choose the next $\zeta=\zeta(k_1,\cdots,k_m)$ so that each $k_i>> k_\gamma$, the 
new $\zeta$ cannot be isomorphic to the ones previously chosen.   

On the other hand, the Seifert fibered components $M_j$ have larger mapping class groups and 
require a little more care.  Let $\pi_j: M_j\rightarrow B_j$ be the Seifert fibration map which 
projects to the base $B_j$.  In case $M_j$ is toroidal --- see the paragraph after Theorem 
\ref{JSJ} for exceptions ---  the torus $T$ is vertical, and there exists another torus $T'$ 
which intersects $T$ {\it persistently}.    This case is already treated in Colin \cite{Co99b}, 
where the torsion for $T$ and $T'$ are shown to be finite.  The two remaining cases are:
$B_j$ is $S^2$ with at most three singular points or punctures, or $B_j$ is $\R\P^2$ with at
most two singular points or punctures.

In the $S^2$ case, if there are two punctures and no singular points, then $M_j$ is $T^2\times
I$ and the minimality of the Torus Decomposition of
$M$ implies that $M$ must be a torus bundle over $S^1$, in which case the  theorem is already
proved by Giroux \cite{Gi99a}.  If there is one puncture and one or no singular
points,  then the torus is compressible.  Therefore, for $S^2$, we assume that there is a
total of three  singular points or punctures (with at least one puncture).  For $\R\P^2$, if
there is one  puncture and no singular points, we have a twisted $I$-bundle over a Klein
bottle which is  separated by the torus boundary.   We will treat this case separately.

Let $S_j$ be a norm-minimizing, oriented, surface on $M_j$ that is not homologous to a vertical
annulus.  It follows that $S_j$ may be made horizontal, and by the Riemann-Hurwitz formula, $S_j$
cannot be an  annulus or a disk. A diffeomorphism $\phi$ of $M_j$ must take $S_j$ to another
horizontal  surface $\phi(S_j)$, since $S_j$ is not an annulus and hence cannot be vertical
(Proposition  \ref{horizontal-vertical}).   A Riemann-Hurwitz  formula calculation reveals 
that $\pi:\phi(S_j)\rightarrow B_j$   and $\pi:S_j\rightarrow B_j$ have the same number of 
sheets.  Moreover, $\phi$ is isotopic to a diffeomorphism which preserves the Seifert 
fibration.  Now, if we construct $\gamma$ as before from adjacent $M_j$ and $M_{j+1}$ 
($M_{j+1}$ not necessarily Seifert fibered), then $t(\mathcal{L}_\gamma)$ remains invariant 
under various choices of $S_j$ up to diffeomorphism.  This is because a diffeomorphism $\phi$ 
of $M_j$ which does not permute the boundary components will induce a map
$T_i=\R^2/\Z^2 \rightarrow T_i= \R^2/\Z^2$ which sends $(x,y)\mapsto (x,px+y)$, $p\in \Z$, on 
each boundary component $T_i$.   The rest of the argument goes through as in the hyperbolic 
case.

Finally assume $M_j$ is a twisted $I$-bundle over a Klein bottle and $T_i=\bdry M_j$.  In this 
case we may distinguish the universally tight contact structures up to isomorphism by passing 
to a double cover $\widetilde M$ of $M$ which is obtained by gluing two copies of $M\backslash 
M_j$ onto the double cover $T^2\times I$ of $M_j$.  (We may need to take larger covers to 
unwind all the twisted $I$-bundle components over Klein bottles.)

The contact structures on $M$ we construct lift to nonisomorphic universally tight contact structures on $\widetilde M$, and therefore it follows that there are infinitely many universally tight contact structures on $M$ up to isomorphism.
\end{proof}

\subsection{Construction of auxiliary manifold $W$}

In this section we construct $W$ satisfying the properties stipulated in the previous section.

\s\n
{\bf Case 1.}
Suppose that $T$ is nonseparating and the next stage of the Gabai 
decomposition comes from
a connected surface $S$ which intersects both $T_+$ and $T_-$.  Let 
$\gamma\subset S$ be an arc
from $T_+$ to $T_-$.   First use residual finiteness to take a large 
finite Galois cover
$\pi:\widetilde M \rightarrow M$ which expands $T$.  Let $\widetilde 
T_1$, $\widetilde T_2$ be
connected components of $\pi^{-1}(T)$, and $\widetilde \gamma\subset 
\widetilde M$ a lift of
$\gamma$ with one endpoint on $\widetilde T_1$ and the other endpoint 
on $\widetilde T_2$.  (By
a ``lift'' $\widetilde\gamma$ of $\gamma$ we mean an arc in $\widetilde 
M$ which maps 1-1 down to
$\gamma$ via $\pi$, with the exception of endpoints of the arc.)  It 
might happen that
$\widetilde T_1=\widetilde T_2$.  In that case, we use residual 
finiteness again to unwrap in
the direction transverse to $\widetilde T_1=\widetilde T_2$, so that 
the two connected
components are distinct.  Now, cut $\widetilde M$ along $\widetilde 
T_1$ and $\widetilde T_2$
and glue together pairs $\widetilde T_{1,+}$, $\widetilde T_{2,-}$ and $\widetilde T_{1,-}$, 
$\widetilde T_{2,+}$ in a $\pi_1(T)$-equivariant manner so that the endpoints of
$\widetilde \gamma$ glue to give a closed curve $\widehat\gamma$ in 
the manifold $W$.   We may need to throw away a component which does not contain 
$\widehat\gamma$.
Note that, in order to glue the boundaries, we need $\widetilde M$ to be a
{\it Galois} cover of $M$ --- i.e., $\pi_1(M)$ must act transitively 
on $\widetilde
M$.  The glued-up contact structure $\zeta$ is universally tight by 
Theorem \ref{variant}.

\s\n
{\bf Case 2.}      Either $T$ is separating and there exist two 
surfaces $S_1$ and $S_2$
or $T$ is nonseparating and there are two surfaces $S_1$ and $S_2$ 
which bound $T_+$, $T_-$ respectively. Let $\gamma_1$, $\gamma_2$ be arcs on $S_1$ and $S_2$ 
which glue to give $\gamma$.  Assume $\gamma_i$ are closed curves (identify their endpoints 
on $T$).  Let $M_i$ be the 
component of $M\backslash T$ containing $S_i$ (note $M_1$ may be the same as $M_2$). 

\begin{lemma}
$\gamma_i\not\in \pi_1(T)$, $i=1,2$, as elements of $\pi_1(M)$.  
\end{lemma}

\begin{proof}
Consider the universal cover 
$\overline\pi:\R^3=\overline M_i\rightarrow M_i$.  We claim that a
lift $\overline\gamma$ of $\gamma=\gamma_i$ has endpoints on different $\R^2$-components of
$\overline\pi^{-1}(T)$. (Call them $\R^2_j$, $j=1,2$.) Assume otherwise; 
namely the  universal cover
$\overline S$ (connected component) of $S=S_i$ intersects a particular 
$\R^2_j$ along at least two
coherently oriented lines.   This follows from the well-grooming of $S$   --- all the 
components of $\bdry S\cap T$ were oriented in the same direction.  Now, $\overline S$ must 
separate $\overline{M_i}=\R^3$ into two half-planes.  However, the coherent 
orientation of $\overline S\cap
\R^2_i$ contradicts the fact that $\overline S$ separates.  This 
implies that $\gamma\not \in
\pi_1(T)\subset \pi_1(M_i)$.   To conclude the proof, we simply note that $\pi_1(M_i)$ injects 
into $\pi_1(M)$. 
\end{proof}

Now we use the fact that $\pi_1(T)$ is separable in $\pi_1(M)$ (Theorem \ref{separability}).     
There exist finite index subgroups $K_1$, $K_2$ in $\pi_1(M)$ which contain $\pi_1(T)$ 
but miss $\gamma_1$, $\gamma_2$, respectively.    The finite cover corresponding to $K_1\cap 
K_2$ will then satisfy the condition that the endpoints of lifts of $\gamma_1$ lie on distinct 
connected torus components of the preimage of $T$, and the same holds for $\gamma_2$.  Finally 
take a finite index normal subgroup $K\subset K_1\cap K_2$ which sufficiently expands $T$.  
Let $\pi:\widetilde M\rightarrow M$ be the corresponding finite cover and $\widetilde 
\gamma=\widetilde \gamma_1+\widetilde\gamma_2$ a lift of $\gamma$.   
Here, let $\widetilde T_j$, $j=1,2,3$, be three components of 
$\pi^{-1}(T)$, and 
$\widetilde \gamma_i$, $i=1,2$, connect from $\widetilde T_i$ to 
$\widetilde T_{i+1}$.  (The endpoint of 
$\widetilde\gamma_1$ equals the initial point of $\widetilde\gamma_2$.)
Now
cut $\widetilde M$ along $\widetilde T_1$ and $\widetilde T_3$ and 
glue pairs $\widetilde T_{1,+}$, $\widetilde T_{3,-}$ and $\widetilde T_{1,-}$, $\widetilde 
T_{3,+}$ equivariantly with
respect to $\pi_1(T)$ so that $\widetilde \gamma$ glues into a closed 
curve $\widehat \gamma$.     This is the desired manifold $W$.

\s
\n
{\it Acknowledgements.}  We thank Emily Hamilton and Darren Long for their help with subgroup 
separability and residual finiteness.    The first author also would like to thank Brian Conrey 
and the American Institute of Mathematics, where part of this paper was written.

\end{document}